\documentclass{article}
\usepackage[utf8]{inputenc}

\usepackage{amsmath}
\usepackage{amssymb}
\usepackage{mathtools}
\usepackage{graphicx}

\title{A Partial Residue Categorization of the Magic Square of Squares}
\author{Chrisitan Woll\\
cwoll@ucsd.edu}
\date{}

\begin{document}

\maketitle

\section{Introduction}

It is unknown at present whether a magic square of squared integers exists. Such an object is defined to be a $3\times3$ grid of $9$ distinct integer squares,
\begin{equation}
\begin{tabular}{c|c|c}
\label{msos}
  $a^2$ & $b^2$ & $c^2$ \\      \hline
  $d^2$ & $e^2$ & $f^2$ \\      \hline
  $g^2$ & $h^2$ & $i^2$
\end{tabular}\ ,
\end{equation}
such that the entries of each row, column, and two main diagonals sum to the same total, say $T$. The problem has a bounty and has attracted the interest of many mathematicians -- mostly recreational [1].\\
\indent A magic square is called \textit{primitive} if the greatest common divisor of all entries is $1$ and is called \textit{reducible} otherwise. Since the entries of any magic square may be scaled to produce another valid magic square, it follows very simply that the existence of any reducible magic square implies the existence of a primitive square (hence the name 'reducible'). In this paper we analyze some properties of the primitive magic square of squares mod $p$ for various primes $p$.

\section{Divisibility of the Central Entry}

An immediate consequence of our definitions is\\\\
\textbf{Theorem 2.1:} The total of any magic square of squares, $T$, is 3 times the central entry, $e^2$.\\\\
\textbf{Proof:} Consider the five sums
$$a^2+e^2+i^2\ =\ b^2+e^2+h^2\ =\ c^2+e^2+g^2\ =\ T,$$
$$\text{and}\quad a^2+b^2+c^2\ =\ g^2+h^2+i^2\ =\ T.$$
From them, we can rearrange the desired proposition,
$$3T=(a^2+e^2+i^2)+(b^2+e^2+h^2)+(c^2+e^2+g^2)$$
$$=(a^2+b^2+c^2)+(e^2+e^2+e^2)+(g^2+h^2+i^2)=2T+3e^2.$$
Subtracting $2T$ yields $T=3e^2$.\\
\indent\hfill$\square$\\\\
Before analyzing residues it is helpful -- though perhaps trivial -- to note\\\\
\textbf{Lemma 2.1:} The central entry in any magic square of squares is divisible by at least one prime.\\\\
\textbf{Proof: }Since the entries are distinct, we may bound the total,
$$3e^2=T\ge 1^2+2^2+3^2=14.$$
It follows that $e>2$. Improving this bound is not difficult.\\
\indent\hfill$\square$

\section{Residues in $\mathbb{F}_p$}

Suppose $a^2,b^2,...,i^2$ to form a valid magic square of squares and $p$ to be a prime dividing $e$. We will now analyze the magic square mod $p$. The symbol $\mathbb{F}_p$ is used to denote the finite field formed by the residues of the integers mod $p$. The overline, $\overline{\ \cdot\ }$, is used to denote the residue of an integer in $\mathbb{F}_p$ (I.e. $\overline{a}=a\ \%\ p$). Thus, the residue of the square in line \ref{msos} is of the form
\begin{equation}
\begin{tabular}{c|c|c}
  $\overline a^2$ & $\overline b^2$ & $\overline c^2$ \\      \hline
  $\overline d^2$ & $0$ & $\overline f^2$ \\      \hline
  $\overline g^2$ & $\overline h^2$ & $\overline i^2$
\end{tabular}\ .
\end{equation}
\indent We call this object the \textit{residue class} of the magic square. Since the additive relations of the square remain valid mod $p$, a residue class may be thought of as a solution to the magic square of squares problem over $\mathbb{F}_p$ instead of over the integers. The residue class lends itself to various derivations. For the remainder of the paper, we let context decide when a number is to be treated as an integer or as a residue in $\mathbb{F}_p$.  \\\\
\textbf{Theorem 3.1:} The central entry of any primitive magic square of squares is divisible only by $p=2$ and primes of the form $p=4k+1$ where $k$ is an integer.\\\\
\textbf{Proof:} It follows from Theorem 2.1 that if $p$ divides $e$, then the square's total, $T$, is divisible by $p$ and thus the total of the residue class, $\overline{T}$, is zero. Consider the sums which pass through the square's center,
\begin{equation}
\label{sos}
    \overline{a}^2+\overline{i}^2\ =\ \overline{b}^2+\overline{h}^2\ =\ \overline{c}^2+\overline{g}^2\ =\ \overline{d}^2+\overline{f}^2\ =\ 0.
\end{equation}
\indent In any primitive magic square, there is no prime dividing all the entries, so at least one of $\overline{a},\overline{b},\overline{c},\overline{g},\overline{h},$ or $\overline{i}$ is nonzero. Therefore we have a nontrivial -- meaning nonzero -- solution to $x^2+y^2=0$ in $\mathbb{F}_p$. It follows that $\big(\frac{x}{y}\big)^2=-1$. Thus, if $-1\not=1$ in $\mathbb{F}_p$, then $\frac{x}{y}$ is an element of order $4$ within $\mathbb{F}_p^\times$, the multiplicative group formed by the nonzero elements of $\mathbb{F}_p$. \\
\indent We continue with an argument from group theory. The element of order $4$ generates a subgroup of order $4$ within the group of nonzero elements of $\mathbb{F}$, denoted here as $\mathbb{F}_p^\times$. This, by Lagrange's Theorem [2], implies $4$ divides the group's order, $|\mathbb{F}_p^\times|=p-1$ and that $p$ is therefore of the form $4k+1$. The case $-1=1$ occurs only in $\mathbb{F}_2$. Thus if $p$ divides the central entry of a primitive magic square of squares, either $p=2$ or $p$ is of the form $4k+1$.\\
\indent Conversely, if $p$ is of the form $4k+3$, $p$ must divide all other entries in the magic square which, by definition, means the square is reducible.\\
\indent\hfill$\square$\\\\
\indent Consider first the case that $p=2$ divides $e$. There are only $4$ possibilities for the residue class of the magic square,
\begin{equation}
\begin{tabular}{c|c|c}
  $0$ & $0$ & $0$ \\      \hline
  $0$ & $0$ & $0$ \\      \hline
  $0$ & $0$ & $0$
\end{tabular}\ ,\quad
\begin{tabular}{c|c|c}
  $1$ & $1$ & $0$ \\      \hline
  $1$ & $0$ & $1$ \\      \hline
  $0$ & $1$ & $1$
\end{tabular}\ ,\quad
\begin{tabular}{c|c|c}
  $1$ & $0$ & $1$ \\      \hline
  $0$ & $0$ & $0$ \\      \hline
  $1$ & $0$ & $1$
\end{tabular}\ ,\quad\text{and}\quad
\begin{tabular}{c|c|c}
  $0$ & $1$ & $1$ \\      \hline
  $1$ & $0$ & $1$ \\      \hline
  $1$ & $1$ & $0$
\end{tabular}\ .
\end{equation}
Interestingly, these $4$ residue classes form a copy of the Klein four-group when added element-wise to each other. \\\\
\textbf{Corollary 3.1:} If the central entry of a magic square is even, then at least one the central column, the central row, or main diagonals will consist of all even numbers.\\\\
\indent It remains to categorize the residue classes for primes of the form $p=4k+1$. In such a case, we let $w$ denote an element -- there are two -- of order $4$ in $\mathbb{F}_p^\times$. The relation $\overline{a}^2+\overline{i}^2=0$ implies $\overline{i}^2=(w\overline{a})^2$. Similarly, $\overline{h}^2=(w\overline{b})^2,\overline{g}^2=(w\overline{c})^2,$ and $\overline{f}^2=(w\overline{d})^2$. Thus the corresponding residue class is parametrized entirely by $\overline{a},\overline{b},\overline{c},$ and $\overline{d}$ and is of the form
\begin{equation}
\begin{tabular}{c|c|c}
\label{4par}
  $\overline a^2$ & $\overline b^2$ & $\overline c^2$ \\      \hline
  $\overline d^2$ & $0$ & $(\overline dw)^2$ \\      \hline
  $(\overline cw)^2$ & $(\overline bw)^2$ & $(\overline aw)^2$
\end{tabular}\ .
\end{equation}
\indent The residue classes for primes of the form $4k+1$ fall nicely into two categories, \textit{trivial} and \textit{nontrivial}.

\section{Trivial Residue Classes}

\indent Suppose $p$ to be a prime of the form $4k+1$ dividing the central entry of a magic square of squares and again let $w$ denote an element of order $4$ in $\mathbb{F}_p$ (meaning $w^2=-1$). The magic square's residue class is called \textit{trivial} if other entries besides the center are zero. \\
\indent If a corner element, say $\overline{a}^2$, is zero, then the residue class is of the form
\begin{equation}
\begin{tabular}{c|c|c}
  $0$ & $\overline b^2$ & $\overline c^2$ \\      \hline
  $\overline d^2$ & $0$ & $(\overline dw)^2$ \\      \hline
  $(\overline cw)^2$ & $(\overline bw)^2$ & $0$
\end{tabular}\ 
\end{equation}
and the remaining relations are $\overline{b}^2+\overline{c}^2=\overline{c}^2+(\overline{d}w)^2=0$. It follows that $\overline{c}^2=\overline{d}^2=(\overline{b}w)^2$. Thus the whole square is parametrized by $\overline{b}$, 
\begin{equation}
\begin{tabular}{c|c|c}
  $0$ & $\overline b^2$ & $(w\overline b)^2$ \\      \hline
  $(w\overline b)^2$ & $0$ & $\overline b^2$ \\      \hline
  $\overline b^2$ & $(w\overline b)^2$ & $0$
\end{tabular}\ .
\end{equation}
\indent Here is a good place for a word about scaling. Just as a reducible magic square of squares can be scaled to a primitive, so too can a residue class be scaled so as to produce another residue class. Consider the former residue class with $\overline{b}=1$ in $\mathbb{F}_{13}$,
\begin{equation}
\begin{tabular}{c|c|c}
  $0$ & $1^2$ & $12^2$ \\      \hline
  $12^2$ & $0$ & $1^2$ \\      \hline
  $1^2$ & $12^2$ & $0$
\end{tabular}\ .
\end{equation}
Each entry may be scaled by a \textit{quadratic residue} which -- to define the term -- means the remainder of any integer square mod $p$ (ex. $4^2\ \%\ 13=3$ means $3$ is a quadratic residue mod $13$). Scaling the former residue class by $3=4^2$ produces
\begin{equation}
\begin{tabular}{c|c|c}
  $0$ & $4^2$ & $6^2$ \\      \hline
  $6^2$ & $0$ & $4^2$ \\      \hline
  $4^2$ & $6^2$ & $0$
\end{tabular}\ .
\end{equation}
\indent Returning to the general case in which a corner is zero, the whole square is parametrized by $\overline{b}$. Since the square is assumed to be primitive, $\overline{b}$ must be nonzero and it follows that the entries of the residue class may be scaled by $1/\overline{b}^2$ to produce,
\begin{equation}
\begin{tabular}{c|c|c}
  $0$ & $1$ & $w^2$ \\      \hline
  $w^2$ & $0$ & $1$ \\      \hline
  $1$ & $w^2$ & $0$
\end{tabular}\ .
\end{equation}\\
\textbf{Corollary 4.1:} If a prime of the form $p=4k+1$ divides both 
the central entry and a corner entry of a magic square of squares, the square's residue class mod $p$ will be of the form
\begin{equation}
\begin{tabular}{c|c|c}
  $0$ & $1$ & $-1$ \\      \hline
  $-1$ & $0$ & $1$ \\      \hline
  $1$ & $-1$ & $0$
\end{tabular}\ 
\end{equation}
up to scaling and rotation.\\\\
\indent Alternatively, if a mid-edge entry, say $\overline{b}^2$, is zero, the remaining relations become $\overline{a}^2+\overline{c}^2=\overline{a}^2+\overline{d}^2-\overline{c}^2=0$. It follows that $2\overline{a}^2+\overline{d}^2=0$ and again that $2=\big(w\overline{d}/\overline{a}\big)^2$. Thus $2$ is a quadratic residue mod $p$ which, it is known, occurs if and only if $p\equiv 1,7\ (\text{mod}\ 8)$ [3]. And since we are already considering only $p\equiv 1\ (\text{mod}\ 4)$, the restrictions mutually becomes $p\equiv 1\ (\text{mod}\ 8)$. \\\\
\textbf{Corollary 4.2:} If a prime of the form $p=8k+5$ divides the central entry of a magic square of squares, then the other entries of the middle row and middle column (i.e. the `mid-edge' entries) are not divisible by $p$.\\\\
\indent So if a prime of the form $4k+1$ divides the central entry and a mid-edge entry in a magic square of squares, then the prime  is actually of the form $8k+1$. Let $\tau$ denote one of the square roots of $2$ in $\mathbb{F}_p$. Then it follows from the aforementioned relations that $\overline{c}^2=(w\overline{a})^2$ and $\overline{d}^2=(w\tau \overline{a})^2$. And thus the residue class is of the form
\begin{equation}
\begin{tabular}{c|c|c}
  $\overline a^2$ & $0$ & $(w\overline a)^2$ \\      \hline
  $(w\tau\overline a)^2$ & $0$ & $(\tau\overline a)^2$ \\      \hline
  $\overline a^2$ & $0$ & $(w\overline a)^2$
\end{tabular}\ .
\end{equation}
\indent Since the square is primitive, $\overline{a}$ must be nonzero and, as before, the residues may be scaled by $1/\overline{a}^2$ producing
\begin{equation}
\begin{tabular}{c|c|c}
  $1$ & $0$ & $w^2$ \\      \hline
  $(w\tau)^2$ & $0$ & $\tau^2$ \\      \hline
  $1$ & $0$ & $w^2$
\end{tabular}\ .
\end{equation}
\textbf{Corollary 4.3:} If a prime divides the central entry and a mid-edge entry of a primitive magic square of squares, $p$ is of the form $8k+1$ and the square's residue class mod $p$ will be of the form
\begin{equation}
\begin{tabular}{c|c|c}
  $1$ & $0$ & $-1$ \\      \hline
  $-2$ & $0$ & $2$ \\      \hline
  $1$ & $0$ & $-1$
\end{tabular}\ 
\end{equation}
up to scaling and rotation.

\section{Nontrivial Residue Classes}

We call a residue class \textit{nontrivial} if only the central entry is zero. The residue class (given on line \ref{4par}) may therefore be scaled by the inverse of $\overline{a}^2,\overline{b}^2,\overline{c}^2,$ or $\overline{d}^2$. We choose $\overline{c}^2$,
\begin{equation}
\begin{tabular}{c|c|c}
  $(\overline a/\overline{c})^2$ & $(\overline b/\overline{c})^2$ & $1$ \\      \hline
  $(\overline d/\overline{c})^2$ & $0$ & $(w\overline d/\overline{c})^2$ \\      \hline
  $w^2$ & $(w\overline b/\overline{c})^2$ & $(w\overline a/\overline{c})^2$
\end{tabular}\ .
\end{equation}
The square may be cleaned up nicely by letting $\alpha=\overline{d}\overline{c},\beta=w\overline{a}/\overline{c},$ and $\gamma=\overline{b}/\overline{c}$,
\begin{equation}
\begin{tabular}{c|c|c}
\label{3par}
  $(w\beta)^2$ & $\gamma^2$ & $1$ \\      \hline
  $\alpha^2$ & $0$ & $(w\alpha)^2$ \\      \hline
  $w^2$ & $(w\gamma)^2$ & $\beta^2$
\end{tabular}\ .
\end{equation}
Thus a residue class may be instantiated by choosing any $\alpha,\beta,\gamma$ in $\mathbb{F}_p$ satisfying $(w\beta)^2+\gamma^2+1=(w\beta)^2+\alpha^2-1=0$, or equivalently, satisfying
\begin{equation}
\label{seq}
    \alpha^2-\beta^2=\beta^2-\gamma^2=1.
\end{equation}\\
\textbf{Corollary 5.1:} There exists a nontrivial residue class mod $p$ if and only if there exists 3 consecutive nonzero quadratic residues mod $p$.\\\\
\indent We use the symbol $S_p$ to denote the quadratic residues mod $p$. This means $S_p$ is the set $\{n^2:n\in\mathbb{F}_p^\times\}$ which is in fact a subgroup of $\mathbb{F}_p^\times$ having order $\frac{1}{2}|\mathbb{F}_p^\times|=\frac{p-1}{2}$. Ex.
$$S_5=\{1^2,2^2,3^2,4^2\}=\{1,4\}$$
\indent It will be easier to continue if the first few $S_p$ are given explicitly (still restricting ourselves to primes of the form $4k+1)$,
$$S_{5}=\{1,4\}$$
$$S_{13}=\{1, 3, 4, 9, 10, 12\}$$
$$S_{17}=\{1,2,4,8,9,13,15,16\}$$
$$S_{29}=\{1,4,5,6,7,9,13,16,20,22,23,24,25,28\}$$
$$S_{37}=\{1,3,4,7,9,10,11,12,16,21,25,26,27,28,30,33,34,36\}$$
\indent Note that $p=29$ is the smallest such prime having $3$ consecutive nonzero residues, $4=2^2,5=11^2,$ and $6=8^2$ (and $7=6^2$ actually). As stated above, the solution $6^2-8^2=8^2-11^2=1$ may be transformed into the nontrivial residue class
\begin{equation}
\begin{tabular}{c|c|c}
  $9^2$ & $11^2$ & $1^2$ \\      \hline
  $6^2$ & $0$ & $14^2$ \\      \hline
  $w^2$ & $16^2$ & $8^2$
\end{tabular}\ .
\end{equation}
\indent Note that there are therefore no nontrivial residue classes for $p=5,13,$ or $17$.\\\\
\textbf{Corollary 5.2:} If the central entry of a primitive magic square of squares is divisible by $2,5,13,$ or $17$, then so is at least one other entry.

\section{Counting Nontrivial Residue Classes}

\indent We saw that nontrivial residue classes exist for $p=29$ but not for $p=5,13,$ or $17$. In this section we will attempt to count how many nontrivial residue classes each prime yields. Let's create a set hold these consecutive triplets,
$$C_p=\{n\in S_p:n+1,n+2\in S_p\}.$$
\indent We have defined $C_p$ such that if $\alpha^2-\beta^2=\beta^2-\gamma^2=1$ for $\alpha,\beta,\gamma\not=0$, then $\gamma^2$ will be a member of $C_p$ and will be a sort of `representative' for it's sequence. We give the first few $C_p$ explicitly,
$$C_2=C_5=C_{17}=\emptyset$$
$$C_{29}=\{4,5,22,23\}$$
$$C_{37}=\{9,10,25,26\}$$
$$C_{41}=\{8,31\}$$
\indent Note that if $n$ is in $C_p$, then so is $-(n+2)$. This is due to the presence of $w$ in $\mathbb{F}_p$ (recalling that $w^2=-1$). From $\alpha^2-\beta^2=\beta^2-\gamma^2=1$, it follows that
\begin{equation}
\label{seqw}
    (w\gamma)^2-(w\beta)^2=(w\beta)^2-(w\alpha)^2=1.
\end{equation}
Thus if $\gamma^2$ is in $C_p$, so is $(w\alpha)^2=-\alpha^2=-(\gamma^2+2)$. Look again at the residue class on line \ref{3par}. Both the sequences, on line \ref{seq} and line \ref{seqw}, occur in the same residue class. Thus this alternative solution we created with scaling-by-$w$, corresponds to a reflection of the magic square about the upwards-slanting main diagonal.\\
\indent A new solution to the equation on line \ref{seq} may also be produced by dividing through by $\beta^2$. This however, corresponds only to the case of scaling the square on line \ref{seq} by $1/\overline{a}^2$ instead of by $1/\overline{c}^2$. There are no interesting consequences.\\\\
\textbf{Theorem 6.1:} The number of residue classes produced by a primitive magic squares of squares whose central entries are divisible by $p$ is no more than
$$(p-1)\cdot\Big(|C_p|+
\begin{cases}
4& p\equiv 1\ (\text{mod}\ 8)\\
2& p\equiv 5\ (\text{mod}\ 8)
\end{cases}\ \Big).
$$\\
\textbf{Proof:} Let $k$ be the number of trivial residue classes for a given prime $p$ up to scaling and rotation. As seen in Section 4 that that $k=2$ for $p\equiv 1\ (\text{mod}\ 8)$ and $k=1$ for $p\equiv 5\ (\text{mod}\ 8)$. Since there are $4$ choices of rotation and $|S_p|=(p-1)/2$ choices for scaling, there are $2k(p-1)$ possible trivial residue classes in total. Note that each of the trivial residue classes is fixed by some reflection.\\
\indent Similarly, there are $\frac{1}{2}|C_p|$ nontrivial residue classes up to scaling and rotation (noting that rotation by $180^\text{o}$ and scaling by $w^2$ fixes the square). Again, there are $4$ choices of rotation and $(p-1)/2$ choices for scaling. Thus there are $4\frac{p-1}{2}\cdot \frac{1}{2}|C_p|=(p-1)|C_p|$ possible nontrivial residue classes.\\
\indent And in total there are no more than $(p-1)(|C_p|+2k)$ residue classes produced by any prime $p$.\\
\indent\hfill$\square$\\\\
\indent The set $C_p$ is nonempty for every prime of the form $4k+1$ with the exception of only $5,13,$ and $17$ [6][7]. The author did not understand the papers proving this to be so. We present here some modular relations which guarantee some $C_p$, but not all, to be nonempty and offer construction of the consecutive residues. We will present one such relation with proof and leave the reader to determine other such relations should they be needed.\\\\
\textbf{Theorem 6.2:} The set $C_p$ is nonempty if $p\equiv 1,9\ (\text{mod}\ 20)$.\\\\
\textbf{Proof:} Consider the integer equation
$$49^2-41^2=41^2-31^2=720.$$
Our goal here is to scale this equation and take the residues mod $p$ to produce a solution in $\mathbb{F}_p$ to the equation on line \ref{seq}. Suppose $p\ge 41$. Then $\overline{49},\overline{41},$ and $\overline{31}$ (the residues in $\mathbb{F}_p$) will be distinct and nonzero. Also, in order to scale the equation properly, $720$ must be a quadratic residue mod $p$. If so then each of the terms will remain squares when scaled by $1/\overline{720}$.\\
\indent Since $720=5\cdot 12^2$, it follows that $720$ is a quadratic residue whenever $5$ is with exceptions $p=2$ and $p=3$ (but we need not worry about these since we have already restricted $p>41$). It is known that $5$ is a quadratic residue exactly when $p\equiv 1,9\ (\text{mod}\ 10)$ [3]. Thus since we are dealing only with the case that $p\equiv 1\ (\text{mod}\ 4)$, the modular relation may be strengthened to $p\equiv 1,9\ (\text{mod}\ 20)$.\\
\indent The two sets exempted, $C_{29}$ and $C_{41}$, can be checked by hand. They were given explicitly above.\\
\indent\hfill$\square$\\\\
\indent Let's get our hands dirty and work through that proof with $p=61$ (meaning the following computations will be done in $\mathbb{F}_p$). \\
\indent The residues of the integer squares are $49^2=22,\ 41^2=34,$ and $31^2=46$. As for the common difference, $5=26^2$. Thus $720=5\cdot12^2=(26\cdot 12)^2=7^2$. So the value by which we scale the whole equation is $1/720=1/7^2$. The inverse of $7$ in $\mathbb{F}_p$ is $35$. Thus $1/7^2=35^2=5$ and -- this last part will be done with written commentary -- we conclude with:
$$49^2-41^2=41^2-31^2=720$$
$$\Rightarrow\quad 22-34=34-46=7^2$$
$$\Rightarrow\quad 22/7^2-34/7^2=34/7^2-46/7^2$$
$$= 22\cdot 5-34\cdot 5=34\cdot 5-26\cdot 5$$
$$= 22\cdot 5-34\cdot 5=34\cdot 5-26\cdot 5$$
$$=49-48=48-49 =1$$\\
\indent The same treatment can be done with any arithmetic sequence of squares. For example, the sequence $7^2-5^2=5^2-1^2=24$ results in \\\\
\textbf{Theorem 6.3:} The set $C_p$ is nonempty if $p\equiv 1,5\ (\text{mod}\ 24)$ with the exception of $p=5$.\\\\
\textbf{Proof:} Left undone.\\\\
\indent These two theorems cover all primes less than $100$ not already decided in this Section. The primes less than $500$ not yet covered are
$$113,137,157,233,257,277,353,373,\ \text{and}\ 397.$$
Quick numerical tests reveal that $C_p$ is nonempty for all these primes [5].\\

$$\textbf{REFERENCES}$$\\\
[1] HTTP://www.multimagie.com/English/SquaresOfSquaresSearch.htm\\\\\
[2] I.N. Herstein, \textit{Abstract Algebra}.\\\\\
[3] Matthew Bates, \textit{Primes of the Form} $x^2+ny^2$. University of Massechusetts, 
\indent Amherst personal homepage.\\\\\
[4] HTTP://en.wikipedia.org/wiki/Congruum \\\\\
[5] SageMath, the Sage Mathematics Software System (Version 6.9),\\
\indent The Sage Developers, 2015, http://www.sagemath.org.\\\\\
[6] Buell \& Hudson, $\textit{On Runs of Consecutive Quadratic Residues}$\\
\indent $\textit{and Quadratic Nonresidues}$.\\\\\
[7] Noam D. Elkies, $\textit{On Runs of Residues}$\\\\\

\end{document}